\documentclass[11pt,a4paper,leqno]{amsart}
\usepackage{amssymb}
\usepackage{url}
\numberwithin{equation}{section}

\newtheoremstyle{theorem}{3pt}{3pt}%
{\it}%         Body font
{}%         Indent amount (empty = no indent, \parindent = para indent)
{\bfseries}% Thm head font  (could be also \sc)
{:}%        Punctuation after thm head
{.5em}%     Space after thm head (\newline = linebreak)
{}%         Thm head spec
\theoremstyle{theorem}
\newtheorem{theorem}{Theorem}[section]

\newtheorem{definition}[theorem]{Definition}

\newtheoremstyle{example}{3pt}{3pt}%
{}%         Body font
{}%         Indent amount (empty = no indent, \parindent = para indent)
{\sc}% Thm head font  (could be also \sc)
{:}%        Punctuation after thm head
{.5em}%     Space after thm head (\newline = linebreak)
{}%         Thm head spec
\theoremstyle{example}

\newtheoremstyle{remark}{3pt}{3pt}%
{}%         Body font
{}%         Indent amount (empty = no indent, \parindent = para indent)
{\sc}% Thm head font  (could be also \sc)
{:}%        Punctuation after thm head
{.5em}%     Space after thm head (\newline = linebreak)
{}%         Thm head spec
\theoremstyle{remark}

\numberwithin{equation}{section}

\newcommand{\hl}{\\\hline}

\newcommand{\thismonth}{\ifcase\month\or
  January\or February\or March\or April\or May\or June\or
  July\or August\or September\or October\or November\or December\fi
  \space\number\year}
\makeatletter
\newcommand{\low}{\@ifnextchar^{}{^{\vphantom x}}}
\newcommand{\high}{\@ifnextchar_{}{_{\vphantom I}}}
\makeatother
\DeclareSymbolFont{script}{U}{eus}{m}{n}
\DeclareSymbolFontAlphabet{\mathscr}{script}
\DeclareMathSymbol{\EuWedge}{0}{script}{"5E}
\DeclareMathAlphabet{\mathrmsl}{OT1}{cmr}{m}{sl}
\newcommand{\rssymb}[2]{\newcommand{#1}{{\mathrmsl{#2}}}}
\newcommand{\calsymb}[2]{\newcommand{#1}{{\mathcal{#2}}}}
\newcommand{\bbsymb}[2]{\newcommand{#1}{{\mathbb{#2}}}}
\newcommand{\lieoper}[2]{\newcommand{#1}{\mathop
  {\mathfrak{#2}\null}\nolimits}}
\newcommand{\oper}[3][n]{\newcommand{#2}{\mathop
  {\mathrm{#3}\null}\ifx n#1\nolimits\else\limits\fi}}
\newcommand{\rsoper}[3][n]{\newcommand{#2}{\mathop
  {\mathrmsl{#3}\null}\ifx n#1\nolimits\else\limits\fi}}
\bbsymb\C{C} \bbsymb\F{F} \bbsymb\HQ{H}\bbsymb\N{N} \bbsymb\Q{Q}
\bbsymb\R{R} \bbsymb\U{U} \bbsymb\V{V} \bbsymb\W{W} \bbsymb\Z{Z}
\bbsymb\bbf{F} \bbsymb\bbk{K} \bbsymb\bbi{I} \bbsymb\bbl{L} \bbsymb\bbo{O}
\bbsymb\bbj{J}
\bbsymb\bby{Y}
\bbsymb\bbp{P}
\bbsymb\bba{A}
\calsymb\cA{A} \calsymb\cB{B} \calsymb\cC{C} \calsymb\cD{D} \calsymb\cE{E}
\calsymb\cF{F} \calsymb\cG{G} \calsymb\cH{H} \calsymb\cI{I} \calsymb\cJ{J}
\calsymb\cK{K} \calsymb\cL{L} \calsymb\cM{M} \calsymb\cN{N} \calsymb\cO{O}
\calsymb\cP{P} \calsymb\cQ{Q} \calsymb\cR{R} \calsymb\cS{S} \calsymb\cT{T}
\calsymb\cU{U} \calsymb\cV{V} \calsymb\cW{W} \calsymb\cX{X} \calsymb\cY{Y}
\calsymb\cZ{Z}
\renewcommand{\geq}{\geqslant} \renewcommand{\leq}{\leqslant}
\oper\End{End} \oper\Hom{Hom}                    % Vector space constructions
\oper\Sym{Sym} \oper\Skew{Skew}
\oper\Aut{Aut}                                   % Group constructions
\oper\GL{GL} \oper\SL{SL}\oper\Symp{Sp}
\oper\CO{CO} \oper\On{O} \oper\SO{SO} \oper\Pin{Pin} \oper\Spin{Spin}
\oper\CU{CU} \oper\Un{U} \oper\SU{SU} \oper\PSU{PSU}
\rsoper\Diff{Diff} \rsoper\SDiff{SDiff}
\lieoper\der{der}                                % Lie algebra constructions
\lieoper\gl{gl} \lieoper\sgl{sl}\lieoper\symp{sp}
\lieoper\co{co} \lieoper\so{so} \lieoper\spin{spin}
\lieoper\cu{cu} \lieoper\un{u}  \lieoper\su{su}
\rsoper\Vect{Vect} \rsoper\Ham{Ham}
\def\la#1{\hbox to #1pc{\leftarrowfill}}
\def\ra#1{\hbox to #1pc{\rightarrowfill}}

\newcommand{\norm}[2][]{|\mkern-2mu|#2|\mkern-2mu|
  _{\lower1pt\hbox{${}_{#1}$}}}
\newcommand{\Norm}[2][]{\bigl|\mkern-3mu\bigr|#2\bigr|\mkern-3mu\bigr|
  _{\lower1pt\hbox{${}_{#1}$}}}

\rsoper\dimn{dim}                           % dimension
\rsoper\grad{grad}                          % gradient
\rsoper\kernel{ker}\rsoper\image{im}        % kernel and image
\rsoper\alt{alt}   \rsoper\sym{sym}         % alternating and symmetric part
\rsoper\Ad{Ad}     \rsoper\ad{ad}           % adjoint action or bundle
\rsoper\CoAd{CoAd} \rsoper\coad{coad}       % coadjoint action
\rsoper\trace{tr}  \rsoper\trfree{tf}       % trace and tracefree part
\rsoper\detm{det}                           % determinant
\rsoper\Vol{Vol}                            % volume
\rsoper\divg{div}                           % divergence
\rsoper\sign{sign}                          % sign function
\rssymb\iden{id}                            % identity
\rssymb\vol{vol}                            % volume element
\oper\Imag{Im}\oper\Real{Re}                % real and imaginary
\newcommand{\sd}{{\raise1pt\hbox{$\scriptscriptstyle +$}}}
\newcommand{\asd}{{\raise1pt\hbox{$\scriptscriptstyle -$}}}
\newcommand{\sdasd}{{\raise1pt\hbox{$\scriptscriptstyle\pm$}}}
\newcommand{\asdsd}{{\raise1pt\hbox{$\scriptscriptstyle\mp$}}}

\rsoper\scal{scal}
\def\kahl/{k\"ahler}
\def\Kahl/{K{\"a}hler}
\begin{document}
\title[Sasaki-Einstein 7-Manifolds, Orlik Polynomials and Homology]
{Sasaki-Einstein 7-Manifolds, Orlik Polynomials and Homology}
\author[R.Gomez]{Ralph R. Gomez}
\address{Swarthmore College}
\email{rgomez1@swarthmore.edu}

\date{\thismonth}
\begin{abstract}
Let $L_f$ be a link of an isolated hypersurface singularity defined by a weighted homogenous polynomial $f.$
In this article, we give ten examples of $2$-connected seven dimensional Sasaki-Einstein manifolds $L_f$ for which
$H_{3}(L_f, \mathbb{Z})$ is completely determined. Using the Boyer-Galicki
construction of links $L_f$ over particular K\"ahler-Einstein orbifolds, we apply a valid
case of Orlik's conjecture to the links $L_f $ so that one is able to explicitly determine $H_{3}(L_f,\mathbb{Z}).$ We give
ten such new examples, all of which have the third Betti number satisfy $10\leq b_{3}(L_{f})\leq 20$.
\end{abstract}
\vspace{2cm}

\maketitle
\vspace{-2mm}

\section{Introduction}
Sasaki-Einstein geometry has enjoyed explosive activity in the last two decades in part due to its relationship to a deep conjecture in string theory called the AdS/CFT conjecture. According to the conjecture, in certain cases, there is a correspondence between Sasaki-Einstein geometry in dimensions five and seven and superconformal field theories in dimension four and three respectively \cite{Sp}. Thus it is worthwhile to exhibit examples of Sasaki-Einstein manifolds particularly in dimension five and seven.
A rich source of constructing Sasaki-Einstein (SE) metrics of positive Ricci curvature pioneered by Boyer and Galicki in \cite{BG2} is via links of isolated hypersurface singularities defined by weighted homogenous polynomials. These smooth manifolds have been used to show the existence of SE metrics on many types of manifolds such as exotic spheres \cite{BGK}, rational homology spheres (\cite{BGN2},\cite{BGN3}) and connected sums of $S^2\times S^3$ \cite{BG2} (see \cite{BG6} for more comprehensive survey.) In general it is very difficult to determine the diffeomorphism or even homeomorphism type of a given link so determining any such geometric or topological data about the link is always helpful. Along these lines, for a given link of dimension $2n-1$, Milnor and Orlik \cite{MilnorOrlick} determined a formula for the $n-1$ Betti number of the link and later on Orlik conjectured a formula \cite{Orlik}(or see section two) for the torsion in $n-1$ integral homology group. This conjecture, called the Orlik conjecture, is known to hold in certain cases. Both of these formulas have been instrumental in extracting some topological data on certain SE manifolds arising as links. For example, based on work of Cheltsov \cite{CHE}, Boyer gave twelve examples \cite{BYR} of SE 7-manifolds arising from links of isolated hypersurface singularities for which the third integral homology group is completely determined. He used Brieskorn-Pham polynomials and Orlik polynomials (see 2.2), both of which are cases in which Orlik's conjecture holds. Inspired by these examples, the main motivation for this article is to find other examples of SE 7-manifolds arising as links generated by Brieskorn-Pham polynomials or Orlik polynomials so that one can explicitly calculate the third integral homology group.\\
\indent Indeed, the main result of the paper is a list of ten examples (see table 3.1) of SE links defined by Orlik polynomials. Because of this, we are then able to calculate the torsion in the third integral homology group explicitly.  In section two, we review the necessary background and in section three we give the table of ten examples together with the third Betti number and explicit forms of $H_3$.

\section{Background}
\indent Define the weighted $\mathbb{C}^{*}$ action on $\mathbb{C}^{n+1}$ defined by
$$(z_{0},...,z_{n})\longmapsto (\lambda^{w_{0}}z_{0},...,\lambda^{w_{n}}z_{n})$$
where $w_{i}$ are the weights which are positive integers and $\lambda \in \mathbb{C}^{*}$. We use the standard notation $\textbf{w}=(w_0,...,w_n)$ to denote a weight vector. In addition, we assume
$$gcd(w_{0},...,w_{n})=1.$$

\begin{definition}
A polynomial $f\in \mathbb{C}[z_{0},...,z_{n}]$ is weighted homogenous if it
satisfies $$f(\lambda^{w_{0}}z_{0},...,\lambda^{w_{n}}z_{n})=\lambda^{d}f(z_{0},...,z_{n})$$
for any $\lambda \in \mathbb{C}^{*}$ and the positive integer $d$ is the degree of $f$.
\end{definition}
The link $L_f$ of an isolated hypersurface singularity defined by a weighted homogenous polynomial $f$ with isolated singularity only at the origina is given by
$$L_{f}=C_{f}\cap S^{2n+1}$$
where $C_f$ is the weighted affine cone defined by $f=0$ in $\mathbb{C}^{n+1}$. By Milnor \cite{Milnor}, $L_{f}$ is a smooth $n-2$ connected manifold of dimension $2n-1$.\\
\indent Recall that a Sasakian structure on a manifold is an odd dimensional Riemannian manifold which admits a normal contact metric structure. (See \cite{BG6} for a comprehensive survey on Sasakian geometry.)
We can now state a key theorem of Boyer and Galicki \cite{BG2} which gives the mechanism for constructing Sasaki-Einstein manifolds via links of isolated hypersurface singularities defined by weighted homogenous polynomials.
\begin{theorem}\cite{BG2}
The link $L_f$ as defined above admits as Sasaki-Einstein structure if and only if the Fano orbifold $\mathcal{Z}_f$ admits a K\"ahler-Einstein orbifold metric of scalar curvature $4n(n+1)$
\end{theorem}
Note that one simply needs to rescale a K\"ahler-Einstein metric of positive scalar curvature to get the desired scalar curvature in the statement of the theorem. We can think of the weighted hypersurface $\mathcal{Z}_{f}$
as the quotient space of the link $L_{f}$ by the locally free circle action where this circle action comes from the weighted Sasakian structure on the link $L_{f}.$ In fact this whole process is summarized in the commutative diagram \cite{BG2}
\begin{equation*}
\begin{matrix}
L_{f} &&\longrightarrow&& S^{2n+1}_{\textbf{w}} \\
 \Big\downarrow \pi &&&&\Big\downarrow\\
\mathcal{Z}_{f} &&\longrightarrow&& \mathbb{P}(\textbf{w}),
\end{matrix}
\end{equation*}
where $S^{2n+1}_{\textbf{w}}$ denotes the unit sphere with a weighted Sasakian structure, $\mathbb{P}(\mathbf{w})$ is weighted projective space coming from the quotient of $S^{2n+1}_{\textbf{w}}$
by a weighted circle action generated from the weighted Sasakian structure. The top horizontal arrow is a Sasakian embedding and the bottom arrow is K\"ahler embedding. Moreover the vertical arrows are orbifold Riemannian submersions.\\
\indent Thus, a mechanism for constructing 2-connected Sasaki-Einstein 7-manifolds  boils down to finding orbifold Fano K\"ahler-Einstein hypersurfaces in weighted projective 4-space $\mathbb{P}(\textbf{w})$. Johnson and Koll\'ar in \cite{JK} construct 4442 Fano orbifolds and of this list, 1936 of these are known to admit orbifold K\"ahler-Einstein metrics. Therefore, by the above construction we state a theorem of Boyer and Galicki:
\begin{theorem}
\cite{BGN2} There exists 1936 2-connected Sasaki-Einstein 7-manifolds realized as links of isolated hypersurface singularities defined by weighted homogenous polynomials.
\end{theorem}
\indent In \cite{BGN2}, the authors were able to determine many from the list of 1936 which yield rational homology 7-spheres and they also determined the order of $H_{3}(L_f, \mathbb{Z})$. In this paper, we identify ten links of isolated hypersurface singularities which can be given by so called Orlik polynomials, thus allowing us to calculate the third integral homology group explicitly. First we need to define some quantities \cite{MilnorOrlick}:
$$u_{i}=\frac{d}{gcd(d,w_{i})}, \hspace{1.5cm} v_{i}=\frac{w_{i}}{gcd(d,w_i)}$$
Let $L_f$ denote a link of an isolated hypersurface singularity defined by a weighted homogenous polynomial. The formula for the Betti number $b_{n-1}(L_f)$ is given by:
$$b_{n-1}(L_{f})=\sum(-1)^{n+1-s}\frac{u_{i_{1}},\ldots u_{i_s}}{v_{i_1}\ldots v_{i_s}\text{lcm}(u_{i_1},\ldots ,u_{i_s})}.$$
Here the sum is over all possible $2^{n+1}$ subsets $\{i_1,\ldots,i_s\}$ of $\{0,\ldots n\}.$\\

\indent For the torsion data, Orlik conjectured \cite{Orlik} that for a given link $L_{f}$ of dimension $2n-1$ one has
\begin{equation}
H_{n-1}(L_{f},\mathbb{Z})_{tor}=\mathbb{Z}_{d_1}\oplus\mathbb{Z}_{d_2}\oplus \cdots \oplus \mathbb{Z}_{d_r}
\end{equation}
We should now review how the $d_{i}$ data are given, using the presentation given in \cite{BG6}.
Given an index set $\{i_1,i_2,....,i_s\}$, define $I$ to be
the set of all of the $2^s$ subsets and let us designate $J$ to be all of the proper subsets.
For each possible subset, we must define (inductively) a pair of numbers $c_{i_{1},...,i_{s}}$
and $k_{i_{1},...,i_{s}}$. For each ordered
subset $\{i_{1},...,i_{s}\}\subset \{0,1,2,...,n\}$ with $i_{1}<i_2<\cdots <i_{s}$
one defines the set of $2^s$ positive integers, beginning with
$c_{\emptyset}=gcd(u_0,...,u_n):$
$$c_{i_{1},...,i_{s}}=\frac{gcd(u_{0},\ldots,\widehat{u}_{i_1},\ldots,\widehat{u}_{i_s},\ldots,u_n)}{\displaystyle \prod_{J}c_{j_{1},\ldots,j_{t}}}.$$
Now, to get the $k's$:
$$k_{i_{1},...,i_{s}}=\epsilon_{n-s+1}\kappa_{i_1,...,i_s}=\epsilon_{n-s+1}\displaystyle\sum_{I}(-1)^{s-t}\frac{u_{j_1}\cdots u_{j_t}}{v_{j_{1}}\cdots v_{j_t}lcm(u_{j_1},\ldots, u_{j_t})}$$
where
$$\epsilon_{n-s+1}=
\begin{cases}
0, & \text{if } n-s+1\text{ is even}\\
1, & \text{if } n-s+1\text{ is odd.}
\end{cases}
$$
Then for each $1\leq j \leq r=\lfloor max\{k_{i_{1},\ldots ,{i_s}}\}\rfloor$ we put
$$d_{j}=\prod_{k_{i_1,\ldots, i_s}\geq j}c_{i_{1},\ldots, i_{s}}.$$

Though the full conjecture is still open 45 years later, it is known to hold in certain cases. If the link is given by either of the polynomials below
\begin{equation}
z_{0}^{a_0}+z_{1}^{a_{1}}+\cdots + z_{n}^{a_n}, \hspace{.5cm} z_{0}^{a_0}+z_{0}z_{1}^{a_{1}}+\cdots + z_{n-1}z_{n}^{a_n}
\end{equation}
then Orlik's conjecture holds (\cite{Orlik}, \cite{OrlikRandall}). The first type of polynomial is called Brieskorn-Pham and the second one is called Orlik. We will discuss these a bit more in the next section.\\
\indent The formulas for the Betti numbers and torsion would indeed be quite tedious to compute by hand, especially when the degree and the weights are large. Fortunately, Evan Thomas developed a program written in C which computes the Betti numbers and the numbers $d_i$, which generates the torsion in $H_{n-1}(L_f,\mathbb{Z})$. Hence if the link is generated by a Brieskorn-Pham polynomial or an Orlik polynomial, then one explicitly knows the torsion in $H_{n-1}$. This program was also used extensively in (\cite{BG5}, \cite{BG6}, \cite{BYR}). I would like to thank Evan Thomas for giving me permission to use the program and to make it available. It is available at my webpage \url{http://blogs.swarthmore.edu/gomez/} under the heading ``Code for homology of links".\\

\section{Examples}
The Johnson and Koll\'ar list available at \url{https://web.math.princeton.edu/~jmjohnso/} under the heading K\"ahler-Einstein and Tiger of Fano orbifolds in weighted projective space $\mathbb{P}(\textbf{w})$ gives the weight vector $\textbf{w}=(w_{0},w_{1},w_{2},w_{3},w_{4})$ with $w_{0}\leq w_{1}\leq w_{2}\leq w_{3}\leq w_{4}$ (which can always be done after an affine change of coordinates) and it indicates if the weighted hypersurface admits an orbifold K\"ahler-Einstein structure. The degree $d$ is given by $d=(w_{0}+\cdots + w_{4})-1$. It is easy to identify whether or not K\"ahler-Einstein orbifolds on the list come from Brieskorn-Pham polynomials since for a given weight vector $\textbf{w}$ on the list, the exponents of the Brieskorn-Pham polynomial would have to be $a_{i}=d/w_i$ for $i=0,...,4$ and therefore one can do a computer search to see if one gets integer results for the exponents. But are there any coming from Orlik polynomials? To get some Orlik examples, one must search among the weighted hypersurfaces in the list of 1936 K\"ahler-Einstein orbifolds and see if the given weights can be represented by Orlik polynomials. This is more difficult than in the Brieskorn-Pham case since the constraints, given in 3.1, are more complicated. The search was done within the range $9\leq w_{0}\leq 11$ where there are 436 Fano orbifolds. Of this lot, 149 Fano orbifolds are known to admit an orbifold Fano K\"ahler-Einstein structure.  Therefore, for a given weight vector $\textbf{w}=(w_{0},w_{1},w_{2},w_{3},w_{4})$ one needs to see if there exists exponents $a_{i}$, in the Orlik polynomials satisfying
\begin{equation}
d=a_{0}w_{0}=w_{0}+w_{1}a_{1}=w_{1}+w_{2}a_{2}=w_{2}+w_{3}a_{3}=w_{3}+w_{4}a_{4}.
\end{equation}
The ten examples were found by hand, checking many different weights against the given conditions. Once they were found, the computer program developed by Evan Thomas was implemented to determine the Betti number and the torsion data. We now give the table of ten examples.

%\indent We now compute the torsion subgroup of $H_{3}(L_{f},\mathbb{Z})$ where $L_f$ is the link defined by the following Orlik polynomial
%$$z_{0}^{11}+z_{0}z_{1}^{75}+z_{1}z_{2}^{5}+z_{2}z_{3}^{2}+z_{3}z_{4}^{2}.$$
%As indicated in the previous section's discussion we need the $\textbf{u}$ and $\textbf{v}$ data which are given by
%$$\textbf{u}=(11,165,825,825,825)\hspace{1cm} \textbf{v}=(1,2,163,331,247).$$
%To obtain the $c$ numbers, we calculate beginning with the $$c_{\varnothing}=11, c_{0}=15, c_{1}=c_{2}=c_{3}=c_{4}=1.$$
%The only other nontrivial $c$ is $c_{0,1}=15$ and all the other $c's$ are equal to $1.$ The associated $k's$ are then given
%by $$k_{\varnothing}=1,k_{0}=0, k_{0,1}=5.$$
%Thus we see that $d_{1}= 55, d_{2}=d_{3}=d_{4}=d_{5}=5.$ Hence, by Orlik's result we have
%$$torH_{3}(L_{4},\mathbb{Z})=\mathbb{Z}_{55}\oplus (\mathbb{Z}_{5})^{4}.$$

\begin{center}
3.2 Table of Orlik Polynomials
\end{center}
\begin{center}\vbox{\[\begin{array}{|c|c|c|c|c|} \hline \textbf{w}=(w_0,w_1,w_2,w_3,w_4) & {\rm
link\hspace{.1cm} L_{f}}
&deg & b_{3} & H_{3}(M,\mathbb{Z}) \hl
(75,10,163,331,247) & z_{0}^{11}+z_{0}z_{1}^{75}+z_{1}z_{2}^{5}+z_{2}z_{3}^{2}+z_{3}z_{4}^2 & 825 & 10 & \mathbb{Z}^{10}\oplus \mathbb{Z}_{55}\oplus (\mathbb{Z}_{5})^4 \hl
(62,124,155,9,85) & z_{0}^{7}+z_{0}z_{1}^{3}+z_{1}z_{2}^2+z_{2}z_{3}^{31}+z_{3}z_{4}^5 & 434 & 12 & \mathbb{Z}^{12}\oplus\mathbb{Z}_{14}\oplus (\mathbb{Z}_{2})^2 \hl
(9,174,467,277,649)& z_{0}^{175}+z_{0}z_{1}^{9}+z_{1}z_{2}^3+z_{2}z_{3}^{4}+z_{3}z_{4}^2   & 1575 & 12& \mathbb{Z}^{12}\oplus\mathbb{Z}_{525}\oplus(\mathbb{Z}_{3})^2\hl
(87,348,145,11,193)& z_{0}^{9}+z_{0}z_{1}^{2}+z_{1}z_{2}^3+z_{2}z_{3}^{58}+z_{3}z_{4}^{4} & 783 & 12 & \mathbb{Z}^{12}\oplus\mathbb{Z}_{27}\oplus \mathbb{Z}_{3}   \hl
(100,350,9,113,229)& z_{0}^{8}+z_{0}z_{1}^{2}+z_{1}z_{2}^{50}+z_{2}z_{3}^{7}+z_{3}z_{4}^{3} & 800 & 14 &\mathbb{Z}^{14}\oplus\mathbb{Z}_{400}  \hl
(9,291,488,181,787)& z_{0}^{195}+z_{0}z_{1}^{6}+z_{1}z_{2}^{3}+z_{2}z_{3}^{7}+z_{3}z_{4}^{2} & 1755 & 14 & \mathbb{Z}^{14}\oplus\mathbb{Z}_{585}\oplus \mathbb{Z}_{3} \hl
(10,164,333,71,253) & z_{0}^{83}+z_{0}z_{1}^{5}+z_{1}z_{2}^{2}+z_{2}z_{3}^{7}+z_{3}z_{4}^{3} & 830& 14 & \mathbb{Z}^{14}\oplus\mathbb{Z}_{166} \hl
(10,540,275,163,103) & z_{0}^{109}+z_{0}z_{1}^{2}+z_{1}z_{2}^{2}+z_{2}z_{3}^{5}+z_{3}z_{4}^{9} & 1090 & 16 & \mathbb{Z}^{16}\oplus\mathbb{Z}_{218}\oplus \mathbb{Z}_{2} \hl
(32,144,11,103,31) & z_{0}^{10}+z_{0}z_{1}^{2}+z_{1}z_{2}^{16}+z_{2}z_{3}^{3}+z_{3}z_{4}^{7} & 320 &18 & \mathbb{Z}^{18}\oplus\mathbb{Z}_{160} \hl
(45,36,27,11,107) & z_{0}^{5}+z_{0}z_{1}^{5}+z_{1}z_{2}^{7}+z_{2}z_{3}^{18}+z_{3}z_{4}^{2} & 225 & 20 & \mathbb{Z}^{20}\oplus\mathbb{Z}_{5}  \hl

\end{array}\]}
\vspace{3mm}
%\parbox{4.00in}{\small Table 1. The bottom entry needs the number theoretic condition}
\end{center}

\section*{Acknowledgements}
I would like to thank Charles Boyer for useful conversations and I would like to express my gratitude to Evan Thomas for allowing me to use the program he developed and for giving me permission to share the code. Part of this article was prepared with the support of a James Michener Fellowship of Swarthmore College.

%\section{A New Regular Sasaki-Einstein Manifold in Dimension Seven}

%There are some difficulties in constructing quasi-regular Sasaki-Einstein
%manifolds. Indeed, there are obstructions to doing so. But, even if
%one does find examples of odd-dimensional Sasaki-Einstein
%manifolds, in general it is quite difficult to determine the
%homeomorphism or even the diffeomorphism type of the manifold.
%In dimension five, under assumptions simple-connectivity one can use
%Smale's classification theorem. But in higher dimensions, there is no
%such theorem.
%\indent However in a paper of Jiang CITE, the author determines the
%homeomorphism types of certain $2$-connected seven manifolds which
%admit regular circle actions.\\

%\indent Consider the following hypersurface:
%$$X_{6}\subset \mathbb{P}(1,1,1,3).$$
%By the work of Cheltsov, this smooth variety is actually K\"ahler-Einstein.
%This variety arises as a Brieskorn-Pham link $z_{0}^6+z_1^6+z_2^6+z_3^6+z_4^2$
%of degree $6$. Since the weights are pairwise relatively prime, by
%proposition 9.3.22 of CITE , $X_6$ is smooth. Hence, $L_{6}\rightarrow X_{6}$
%is a regular circle bundle, since it admits a regular Sasakian structure.
%Furthermore, by CITE, we have that
%$$H_{3}(L_6,\mathbb{Z})=\mathbb{Z}^{104}\oplus \mathbb{Z}_2.$$
%In particular, this the rank of $H_{3}$ is even and it has torsion. Moreover,
%by Milnor CITE, $L_{6}$ is $2-$connected.

\end{document}